\documentclass{amsart}
\usepackage{latexsym}
\usepackage{amssymb}
\usepackage{amsmath}
\usepackage{amscd}
\usepackage{amsthm}
\usepackage[utf8]{inputenc}
\usepackage{placeins}
\usepackage{fourier}
\usepackage[mathscr]{euscript}
\usepackage{tikz}
\usetikzlibrary{arrows,decorations.markings,%
                shapes,positioning}
                
\newcommand{\ncmd}{\newcommand}
\ncmd{\btk}{\begin{center}\begin{tikzpicture}}
\ncmd{\etk}{\end{tikzpicture}\end{center}}
\ncmd{\dsum}{\displaystyle\sum}
\ncmd{\dprod}{\displaystyle\prod}
\ncmd{\ovl}{\overline}
\ncmd{\ptl}{\partial}
\ncmd{\rar}{\rightarrow}
\ncmd{\mA}{\script{A}}
\ncmd{\mC}{\mathbb{C}}
\ncmd{\mg}{\mathfrak{g}}
\ncmd{\mh}{\mathfrak{h}}
\ncmd{\mH}{\mathbb{H}}
\ncmd{\mk}{\mathfrak{k}}
\ncmd{\mO}{\mathbb{O}}
\ncmd{\mR}{\mathbb{R}}
\ncmd{\mZ}{\mathbb{Z}}
\ncmd{\wdh}{\widehat}
\newtheorem{theorem}{Theorem}
\newtheorem{conjecture}{Conjecture}
\tikzset{E1/.style = {thick, double = black, double distance = 1pt}}
\tikzset{EA/.style = {->,thick, double = black, double distance = 1pt}}
\tikzset{E2/.style = {double = black, double distance = 0.5pt}}
\tikzset{ES/.style = {dashed,double = black, double distance = 0.5pt}}
\tikzset{VS/.style = {circle, inner sep = 2pt, fill = white, very thick, draw = black}}
\tikzset{V1/.style = {circle, inner sep = 1pt, fill = black}}
\title{Invariant Adjoint Tensors of the Classical Groups}
\author{Matthew A. Tai}
\date{}

\begin{document}

\maketitle
\section{Abstract}
For $\mg$ a simple Lie algebra and $G$ its adjoint group, the Chevalley map and work of Coxeter gives a concrete description of the algebra of $G$-invariant polynomials on $\mg$ in terms of traces over various representations. Here we provide an extension of this description to $G$-invariant tensors on $\mg$, although restricted to only providing generators and only for the classical Lie algebras.
\section{Introduction}
For a group $G$ and a $\mC$-linear representation $V$ of $G$, we can look at the polynomial functions on $V$, here denoted $\mC(V)$, and then at the $G$-invariant polynomials $\mC(V)^G$. Chevalley showed that for a simple Lie algebra $\mg$ with adjoint group $G$, Cartan subalgebra $\mh$ and Weyl group $W$, we have that
\begin{equation}
\mC(\mg)^G \cong \mC(\mh)^W
\end{equation}
This map was extended by Harish-Chandra to a relation on the universal enveloping algebra of $\mg$, and was used by Kostant to provide part of his theorem on the decomposition of $\mC(\mg)$ into $G$-representations [Ko63].\\
By work of Coxeter, we know that $\mC(\mh)^W$ is a finitely-generated polynomial algebra with with a minimal generating set of $r = \dim(\mh) = \text{rk}(\mg)$ generators. The corresponding generators of $\mC(\mg)^G$ can be easily described for all $\mg$. For a representation $(V,\pi)$ of a simple Lie algebra $\mg$, we define
\begin{equation}M_V:= \pi(X_\alpha)\otimes K^{\alpha\beta}X_\beta \in End(V)\otimes \mC(\mg)\end{equation}
where the $X_\alpha$ provide a basis for $\mg$ and $K$ is the Killing form in this basis. Then $\mC(\mg)^G$ can be generated by $r$ elements of the form
\begin{equation}tr_V(M_V^k) \in \mC(\mg)\end{equation}
for various $V$ and $k$. The representations $V$ can be almost arbitrary but the set of numbers $k-1$ are fixed by $G$ and are called the exponents $e_i$ of $G$.\\
Here we provide a similar result for the free algebra on $\mg$, which we view as the tensor algebra $T(\mg) = \bigoplus_{k=0}^\infty \mg^{\otimes k}$. We show that the subalgebra $T(\mg)^G$ is generated by elements of a certain form. The set of elements we give is neither finite nor minimal, and we do not provide enough relations to furnish a presentation.
\section{Main Theorem}
For each classical Lie algebra we take the defining representation $(V,\pi)$. For $A_r$, this is one of the two $r+1$-dimensional representations; for $B_r$ this is the $2r+1$-dimensional representation; for $C_r$ and $D_r$ this is the $2r$-dimensional representation. We assume each representation comes with a basis so that $\pi$ yields matrices of the appropriate dimension.\\
Given a degree $k$ tensor $T$ in $T(\mg)^G$, we get an action of the symmetric group $\mathfrak{S}_k$ on $T$ via
\begin{equation}(\sigma.T)(X_1,\ldots,X_k) = T(X_{\sigma^{-1}(1)},\ldots,X_{\sigma^{-1}(k)})\end{equation}
We call this action a permutation of the indices.\\
We define a trace to be a tensor of the form
\begin{equation}T_{V,k}(X_1,\ldots,X_k) = tr_V(\pi(X_1)\pi(X_2)\cdots\pi(X_k))\end{equation}
Define $[-]^a_b:End(V)\rar \mC$ to take an element of $End(V)$ to the $(a,b)$ entry of the matrix corresponding to that element. We define $\pi^k$ to be a degree-$k$ $End(V)$-valued tensor that such that $[\pi^k(-)]_b^a$ takes $(X_1,X_2,\ldots,X_k)$ to 
\begin{equation}[\pi(X_1)]_b^{c_1}[\pi(X_2)]_{c_1}^{c_2}\cdots[\pi(X_k)]_{c_{k-1}}^a\end{equation}
Hence $T_{V,k} = tr_V(\pi^k)$.
\begin{theorem}
For $A_r, B_r,$ and $C_r$, $T(\mg)$ is generated as a tensor algebra by traces $T_V$, with $V$ being the defining representation listed above, allowing for permutation of the indices of tensors. For $D_r$, $T(\mg)$ is generated as tensor algebra by traces $T_V$, with $V$ being the defining representation listed above, and by tensors of the form
\begin{equation}\epsilon_{a_1,\ldots,a_{2r}}\prod_{i=1}^r g^{a_i,b_i}[\pi^{k_i}]_{b_i}^{a_{r+i}}\end{equation}
allowing for permutation of the indices, where $\epsilon_{a_1,\ldots,a_{2r}}$ is the Levi-Civita tensor.
\end{theorem}
Note that in the $D_r$ case, for $k_1 = k_2 = \ldots = k_r = 1$, the extra generating term becomes the Pfaffian up to some phase convention.\\
Remark: For $A_r, C_r$ and $D_r$, the representations chosen here are not actually representations of the adjoint group $G$. Rather they are representations of some cover of $G$ in terms of which all representations of $G$ can be expressed. Given a cover $\tilde{G}$ of $G$, every representation of $G$ lifts to a representation of $\tilde{G}$, so we can pick a cover $\tilde{G}$ that is easier to work with than $G$ and consider all representations of $\tilde{G}$. For $B_r$, the adjoint group is $SO(2r+1,\mC)$, so we can use the $2r+1$-dimensional representation and discuss invariants of $G$ directly.
\section{Proof of the main theorem}
Elements of $T(\mg)^G$ are all formal sums of elements of the form
\begin{equation}N_{\alpha_1,\ldots,\alpha_n} = P_{b_1,\ldots,b_n}^{a_1,\ldots,a_n}[\pi(X_{\alpha_1})]_{a_1}^{b_1}\cdots[\pi(X_{\alpha_n})]_{a_n}^{b_n}\end{equation}
where the $a_i$ and $b_i$ are indices in some representations of $G$ and $P$ is an invariant tensor. Thus the problem comes down to determining what $P$ could be. These involve the invariants on $G$-invariant tensors on $V$. For a proof that the invariants mentioned are the only ones, see [GW09].
\subsection{$A_r$}
For $A_r$, the adjoint group is $PSL(r+1,\mC)$, and all of the representations of $PSL(r+1,\mC)$ can be written as representations of $SL(r+1,\mC)$. In turn, these representations can all be written as symmetric powers or alternating powers or some combination thereof of the $r+1$-dimensional representation $V$, so we only need to consider tensors written in terms of $V$. Thus we only need to consider $a_i$ and $b_i$ to be in $V$ or its dual $V^\vee$. We interpret all of the upper indices as being in $V$ and the lower ones as being in $V^\vee$, with the indices $\alpha_i$ being refered to only as adjoint indices.\\
$V$ contains as invariants only the Kronecker delta, the Levi-Civita tensor and its dual. Hence any invariant tensor $P$ must be built out of Kronecker deltas and Levi-Civita tensors. Each instance of the Levi-Civita tensor introduces $r+1$ lower indices; since there cannot be any free lower $V$ indices in an element of $T(\mg)^G$, each Levi-Civita tensor must be matched by an appropriate source of upper indices, i.e. a dual Levi-Civita tensor.\\
Then we have the identity:
\begin{equation}\epsilon_{a_1,\ldots,a_{r+1}}\epsilon^{b_1,\ldots,b_{r+1}} = \delta_{[a_1}^{b_1}\delta_{a_2}^{b_2}\cdots\delta_{a_{r+1}]}^{b_{r+1}}\end{equation}
which allows us to replace a Levi-Civita and dual Levi-Civita pair with Kronecker deltas. Hence we can reduce any invariant $P$ to a sum of products of Kronecker deltas.\\
Thus
\begin{equation}N_{\alpha_1,\ldots,\alpha_n} = \dsum_{\sigma \in S_n}c_\sigma \delta_{a_{\sigma(1)}}^{b_1}\ldots\delta_{a_{\sigma(n)}}^{b_n}[\pi(X_{\alpha_1})]_{a_1}^{b_1}\cdots[\pi(X_{\alpha_n})]_{a_n}^{b_n}\end{equation}
for some set of coefficients $c_\sigma$. For a given $\sigma$, the contraction of Kronecker deltas and $[\pi(X_{\alpha_i})]$ yields a trace, as desired.
\subsection{$B_r$}
For $B_r$, the adjoint group is $SO(2r+1,\mC)$, and the representations of $SO(2r+1,\mC)$ can be written in terms of the $2r+1$-dimensional representation $V$. Note that while $B_r$ has spin representations which cannot be written as symmetric or alternating powers of $V$, $SO(2r+1,\mC)$ does not have such representations; the spin representations are only representations of $Spin(2r+1,\mC)$.\\
By the same argument as for $A_r$, the invariant $P$ for $B_r$ can be written in terms of the invariants of $V$. $V$ has the Kronecker delta and the Levi-Civita tensors as $A_r$ does, but also comes with the metric $g_{a,b}$ and dual metric $g^{a,b}$. Of these, we only want the Kronecker delta to be necessary for $P$.\\
Unlike in the $A_r$ case, the dual Levi-Civita tensor is not the only source of upper indices that don't come with lower indices. We can use the dual metric $g^{a,b}$ for that. But $g^{a,b}$ gives pairs of upper indices, while each Levi-Civita tensor comes with an odd number of lower indices, so not all of the lower indices from a Levi-Civita tensor can be contracted to indices from dual metrics $g^{a,b}$. Hence there must be at least one other Levi-Civita tensor, which we can turn into a dual Levi-Civita tensor by contracting with copies of $g^{a,b}$, and then the Levi-Civita and dual Levi-Civita pair can be replaced by Kronecker deltas, as in the $A_r$ case.\\
Hence we only have Kronecker deltas and metric/dual metric tensors.\\
Each copy of the metric tensor gives us two lower indices that must be matched by upper indices coming from a dual metric tensor. Unlike in the case of Levi-Civita tensors, the only way to take a metric tensor and a dual metric tensor and replace them by Kronecker deltas is to contract them:
\begin{equation}g_{a,b}g^{b,c} = \delta_a^c = g^{c,b}g_{b,a}\end{equation}
But while $P$ must have a dual metric tensor for each metric tensor it has, the metric tensor might be contracted to a Lie algebra element rather than directly to the dual metric tensor. Consider the following:
\begin{equation}g_{a,b}[\pi(X_1)]_{c_1}^b[\pi(X_2)]_{c_2}^{c_1}\cdots[\pi(X_k)]_{c}^{c_{k-1}}g^{c,d}\cdots\end{equation}
where we know that $g^{c,d}$ must appear somewhere in this expression as otherwise we would have more lower indices than upper indices.\\
The statement that $g_{a,b}$ is invariant under the action of $B_r$ is equivalent to the following:
\begin{equation}[\pi(X)]_a^bg_{b,c} = -g_{a,b}[\pi(X)]_c^b\end{equation}
So we can swap the metric and a Lie algebra element at the cost of a sign change. Hence we can move the metric along the chain until it is contracted directly to a copy of the dual metric, and then replace both by a Kronecker delta. Since every metric and dual metric pair can be eliminated this way, we are left with only Kronecker deltas, which again leaves us with the trace.
\subsection{$C_r$}
For $C_r$, the adjoint group is $PSp(r,\mC)$, but the representations of $PSp(r,\mC)$ can be written as symmetric and antisymmetric powers of the $2r$-dimensional representation of $Sp(r,\mC)$, so we shall do so. The $2r$-dimensional representation has as invariants the Kronecker delta and the Levi-Civita tensors, as well as the symplectic form $f_{a,b}$ and its dual $f^{a,b}$. In fact, the Levi-Civita tensor can be written in terms of the symplectic form by taking $r$ copies of the symplectic form and antisymmetrizing over all of the indices. This is why there is no "special symplectic group"; keeping the symplectic form invariant automatically makes the determinant invariant.\\
So we only need to eliminate the symplectic form. The symplectic form obeys the same invariant equation as the metric did for $B_r$:
\begin{equation}[\pi(X)]_a^bf_{b,c} = -f_{a,b}[\pi(X)]_c^b\end{equation}
so we can swap the symplectic and Lie algebra elements at the cost of a sign change. Hence we can move symplectic forms along the chain until they are contracted to dual symplectic forms, and then convert such pairs to Kronecker deltas. So again we get traces.
\subsection{$D_r$}
The $D_r$ case is somewhat more involved. The adjoint group of $D_r$ is $PSO(2r,\mC)$; unlike the $B_r$ case, $SO(2r,\mC)$ has nontrivial center. The representations of $PSO(2r,\mC)$ can still be written in terms of the $2r$-dimensional representation of $SO(2r,\mC)$, and thus we write everything in terms of the $2r$-dimensional representation. The invariants are the Kronecker deltas, Levi-Civita tensors, and metric tensors, as in the $B_r$ case. Here, though, we cannot rule out the need for Levi-Civita tensors, because here each Levi-Civita tensor gives $2r$ lower indices, and hence the lower indices can be matched entirely by upper indices from the dual metric. However pairs of Levi-Civita tensors can still be replaced by Kronecker deltas, and chains involving only the metric and dual metric but not the Levi-Civita tensors can also be turned into traces.\\
So we have traces as in the $A_r, B_r$ and $C_r$ cases, as well as invariants built from single Levi-Civita tensors and metric tensors.\\
For an invariant involving a single instance of the Levi-Civita tensor, since there are no free lower indices we must again get chains leading from an index on the Levi-Civita tensor, through several Lie algebra elements and copies of the metric and dual metric tensors, and then back to the Levi-Civita tensor by a different lower index. On each chain there must be one more dual metric tensor than metric tensor since the Levi-Civita tensor provides two lower indices. We can use the invariance equation for $g_{a,b}$ to move it past Lie algebra elements until it is contracted with a dual metric tensor and then replace both by a Kronecker delta. So we are left with chains leading from an index on the Levi-Civita tensor, through several Lie algebra elements, and then through exactly one dual metric tensor before connecting to the Levi-Civita tensor again, as described in the theorem.\\
Finally, we need to show that the elements listed in the theorem are not just in $T(\mg)^{\tilde{G}}$ but in $T(\mg)^G$ proper. Here we simply note that the action of $\tilde{G}$ on $\mg$ factors through the action of $G$, and hence any element of $T(\mg)$ that is invariant under the action of $\tilde{G}$ is necessarily invariant under the action of $G$. Hence all of the elements listed in the theorem are in $T(\mg)^G$, despite being written in terms of representations of $\tilde{G}$ rather than representations of $G$.
While the proof is done in terms of representations of $\tilde{G}$, we can also consider traces over representations of $G$. Here we have the issue that the invariants of representations of $G$ are of higher order, cubic or quartic or higher, and do not have the nice relation that the Levi-Civita tensor does. For instance, the adjoint representation has the Kronecker delta, a symmetric bilinear form from the Killing form, the Levi-Civita tensor, and at least one other invariant in the form of the structure constants $f_{\alpha,\beta}^\gamma$, from $[x_\alpha,x_\beta] = f_{\alpha,\beta}^\gamma$. The structure constants obey the Jacobi identity, but the Jacobi identity only involves a specific combination of structure constants, unlike for example the Levi-Civita identity which eliminates a pair regardless of how they are related.\\
For $D_{2k}$, the use of an additional generator type is necessary regardless of what representation is used. Taking a trace and symmetrizing yields an element of $\mC(\mg)$, which by Kostant's result can be written as a polynomial in some set of generators. Given a tensor that is a tensor product of multiple traces, symmetrizing yields the product of the symmetrization of the individual traces. Hence for each generator we only have to consider elements obtained by symmetrizing a single trace rather than symmetrizing tensor products of traces. For a single trace, there is only one symmetrization, so we can distinguish these symmetrizations of single traces by degree. However for $D_{2k}$, the polynomial algebra has two generators of degree $2k$, not one, and so we would need two traces of degree $2k$ with distinct symmetrizations, which cannot occur.\\
For $r = 2k+1$, one could consider one of the spin representations of $D_r$; at least in terms of symmetrization, traces in the spin representation do yield all of the generators of the symmetric algebra, since the Pfaffian in $D_{2k+1}$ is degree $2k+1$ and thus does not overlap with any of the other generators, which all have even degree. However the spin representations yield higher-order invariants, and thus the argument used above does not necessarily hold. Hence we stay with the current description in terms of the $2r$-dimensional representation and with the additional generators.
\section{$GL(n,\mC)$ and $O(n,\mC)$}
The groups $GL(n,\mC)$ and $O(n,\mC)$ also act on $sl(n,\mC)$ and $so(n,\mC)$. $GL(n,\mC)$ and $O(2r+1,\mC)$ are just direct products of $SL(n,\mC)$ and $SO(2r+1,\mC)$ with $\mC^\times$ and $\mZ_2$, so the conjugation action doesn't change.\\
For the case of $SO(2r,\mC)$, $O(2r,\mC)$ is not a central extension, so the conjugation action by elements of the other component does change some behavior. It is still an automorphism of the Lie algebra, though.\\
Due to being an automorphism of the Lie algebra, we get that while $O(2r,\mC)$ doesn't preserve the Pfaffian, which is linear in each basis vector and thus gets a sign change when we change the sign of a single basis vector, the adjoint tensors defined using the Pfaffian are still invariant under $O(2r,\mC)$, as all of the indices in $V$ are contracted over. Hence we get that the adjoint tensor invariants of $GL(n,\mC)$ and $O(n,\mC)$ are identical to those of $SL(n,\mC)$ and $SO(n,\mC)$; despite not being an invariant of $O(2r,\mC)$, the Pfaffian is still an invariant of the Lie algebra $so(2r,\mC)$ and hence gives adjoint tensor invariants.
\section{Some relations}
The Jacobi Identity for the structure constants provide a relation for the traces as listed. In particular, for each Lie algebra the structure constant can be written as the antisymmetrization of a degree $3$ trace and the contraction of two structure constants can be written as a combination of traces in degrees $2$ and $4$. Thus the Jacobi identity, which states that a linear combination of three contractions of pairs of structure constants vanishes, can be written as the vanishing of a linear combination of traces of degree 4 and tensor products of traces in degree 2. The exact statement varies based on $\mg$.\\
Other relations come from the absence of generators of the symmetric algebra in various degrees. If a trace invariant $N$ of $\mg$ has degree $k$ and the polynomial algebra $\mC(\mg)$ has no generator in degree $k$, then the symmetrization of $N$ can be written as a product of symmetrizations of tensor invariants of lower order. For instance for $\mg$ being $B_r$ and $C_r$, there are no elements of $\mC(\mg)$ with odd degree, so the symmetrizations of odd degree traces all vanish identically.\\
Another source of relations is what one might call the Cayley-Hamilton identity for $M$. $M$ is a $\dim(V)$-dimensional square matrix, and hence if we interpret $M$ as living in $End(V)\otimes \mC(\mg)$ rather than in $End(V)\otimes T(\mg)$ then $M$ obeys the Cayley-Hamilton identity, and thus $M^{\dim(V)}$ can be written as a linear combination of $P_kM^k$ for $k < \dim(V)$, where $P_k$ is a polynomial in symmetrized traces of powers of $M$. Thus symmetrized instances of $M^{\dim(V)}$ can be replaced by lower powers of $M$ times other traces. See [RSV] for details of how symmetric traces of other degrees decompose.\\
The author does not claim that all of the relations can be generated from the ones mentioned here.
\section{Diagrams}
We can write the main theorem diagrammatically. Consider the form $V\otimes V^\vee\otimes \mg \rar \mC$ given by the adjoint of the representation $\pi: \mg \rar V\otimes V^\vee$. We denote this form by an triangle, with the thin line corresponding to the adjoint index, the point toward the index in $V^\vee$ and the flat side toward the index in $V$:
\btk
\draw[E2] (0,0)to(0,-1);
\draw[E1] (1, 0) to (-0.1,0);
\draw[E1] (0,0)to (-1,0);
\filldraw[fill = white] (-0.14,0)to(0.1,0.1)to(0.1,-0.1)to(-0.14,0);
\etk
In this notation, a contraction of a $V$ index with a $V^\vee$ index corresponds to joining the ends of the corresponding lines, so that a product of two elements in $\pi(\mg)$ as elements of $End(V)$ is written as:
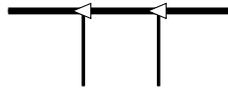
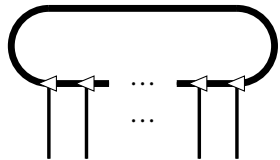
\btk
\draw[E2] (0,0)to(0,-1);
\draw[E2] (-1,0)to(-1,-1);
\draw[E1] (1,0)to (-2,0);
\filldraw[fill = white] (-0.14,0)to(0.1,0.1)to(0.1,-0.1)to(-0.14,0);
\filldraw[fill = white] (-1.14,0)to(-0.9,0.1)to(-0.9,-0.1)to(-1.14,0);
\etk
A trace then becomes
\btk
\draw[E2] (0,0)to(0,-1);
\draw[E2] (-0.5,0)to(-0.5,-1);
\draw[E2] (-2,0)to(-2,-1);
\draw[E2] (-2.5,0)to(-2.5,-1);
\draw[E1] (0,0)to (-0.8,0);
\draw[E1] (-1.7,0) to (-2.5,0);
\draw[E1] (0,0) arc[radius = .5, start angle = -90, end angle = 90];
\draw[E1] (0,1) to (-2.5,1);
\draw[E1] (-2.5,0) arc[radius = 0.5, start angle = 270, end angle = 90];
\node at(-1.25,0){$\ldots$};
\node at(-1.25,-0.5){$\ldots$};
\filldraw[fill = white] (-0.14,0)to(0.1,0.1)to(0.1,-0.1)to(-0.14,0);
\filldraw[fill = white] (-0.5-0.14,0)to(-0.5+0.1,0.1)to(-0.5+0.1,-0.1)to(-0.5-0.14,0);
\filldraw[fill = white] (-2.14,0)to(-1.9,0.1)to(-1.9,-0.1)to(-2.14,0);
\filldraw[fill = white] (-0.5-2.14,0)to(-0.5-1.9,0.1)to(-0.5-1.9,-0.1)to(-0.5-2.14,0);
\etk
and the tensor product of two traces is then two such diagrams placed side by side.\\
The operation of permuting the adjoint indices is then rearranging the free ends of the thin lines. The main theorem then becomes the statement that for $\mg = A_r, B_r,C_r$, any element of $(T(\mg))^G$ can be written as a formal sum of such diagrams.\\
Note that the process of simplifying tensors to traces, we could end up with the following picture:
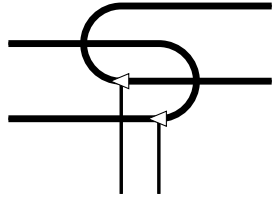
\btk
\draw[E1] (-1,0)to(1,0) arc[radius=0.5,start angle = 90, end angle = -90];
\draw[E1](1,-1)to(-1,-1);
\draw[E1] (2.5,0.5)to(0.5,0.5) arc[radius=0.5,start angle = 90, end angle = 270];
\draw[E1] (0.5,-0.5) to (2.5,-0.5);
\draw[E2] (0.5,-0.5) to (0.5,-2);
\draw[E2] (1,-1)to (1,-2);
\filldraw[fill = white] (1-0.14,-1)to(1.1,-1+0.1)to(1.1,-1-0.1)to(1-0.14,-1);
\filldraw[fill = white] (0.5-0.14,-0.5)to(0.5+.1,-0.5+0.1)to(0.5+.1,-0.5-0.1)to(0.5-0.14,-0.5);
\etk
which we can rewrite as
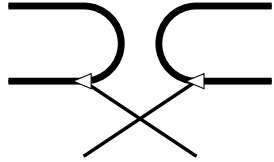
\btk
\draw[E1] (-1,0)to(0,0) arc[radius=0.5,start angle = 90, end angle = -90];
\draw[E1](0,-1)to(-1,-1);
\draw[E1] (2.5,0)to(1.5,0) arc[radius = 0.5, start angle = 90, end angle = 270];
\draw[E1](1.5,-1)to(2.5,-1);
\draw[E2](0,-1)to(1.5,-2);
\draw[E2](1.5,-1)to(0,-2);
\filldraw[fill = white] (0-0.14,-1)to(0+.1,-1+0.1)to(0+.1,-1-0.1)to(0-0.14,-1);
\filldraw[fill = white] (1.5-0.14,-1)to(1.5+.1,-1+0.1)to(1.5+.1,-1-0.1)to(1.5-0.14,-1);
\etk
So any crossings can be transferred to the thin lines, i.e. can be turned into permutation of the adjoint indices.\\
For the case of $D_r$, we use the following notation for the Levi-Civita tensor:
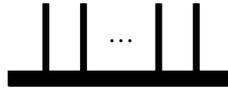
\btk
\filldraw[fill = black] (0,-0.1)to(0,0.1)to(3,0.1)to(3,-0.1);
\draw[E1](0.5,0)to(0.5,1);
\draw[E1](1.0,0)to(1.0,1);
\node at(1.5,.5){$\ldots$};
\draw[E1](2.0,0)to(2.0,1);
\draw[E1](2.5,0)to(2.5,1);
\etk
in accordance to [CV08]. The corresponding adjoint tensors are of the form:
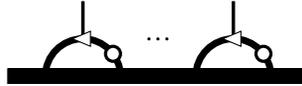
\btk
\filldraw[fill = black] (0,-0.1)to(0,0.1)to(4,0.1)to(4,-0.1);
\draw[E1] (1.5,0) arc[radius= .5,start angle = 0, end angle = 180];
\node[VS]at(1.4,.3){};
\draw[E2](1,0.5)to(1,1);
\filldraw[fill = white] (1-0.14,0.5)to(1+.1,0.5+0.1)to(1+.1,0.5-0.1)to(1-0.14,0.5);
\node at(2.0,.5){$\ldots$};
\draw[E1] (3.5,0) arc[radius= .5,start angle = 0, end angle = 180];
\node[VS]at(3.4,.3){};
\draw[E2](3,0.5)to(3,1);
\filldraw[fill = white] (3-0.14,0.5)to(3+.1,0.5+0.1)to(3+.1,0.5-0.1)to(3-0.14,0.5);
\etk
where the white circles indicate the metric and the thin lines are shorthand for multiple adjoint lines.
\section{Trees}
Note that, up to a scaling, the Killing form on $\mg$ can be written as a trace:
\begin{equation}K_{\alpha\beta} = K(X_\alpha,X_\beta) = tr_V(\pi(X_\alpha)\pi(X_\beta))\end{equation}
The rule for $\pi$ being a representation is that
\begin{equation}\pi(X)\pi(Y) - \pi(Y)\pi(X) = \pi([X,Y])\end{equation}
This allows us to write the structure constants in terms of traces as well:
\begin{eqnarray}c_{\alpha\beta}^\gamma &=& K([X_\alpha,X_\beta],X_\delta)K^{\delta\gamma}\\
&=& tr_V(\pi(X_\alpha)\pi(X_\beta)\pi(X_\delta)-\pi(X_\beta)\pi(X_\alpha)\pi(X_\delta))K^{\delta\gamma}
\end{eqnarray}
It also allows us to write 
\begin{equation}
\pi(X_\alpha)\pi(X_\beta) = \frac{1}{2}(\pi(X_\alpha)\pi(X_\beta)+\pi(X_\beta)\pi(X_\alpha)) + \frac{1}{2}c_{\alpha\beta}^\gamma \pi(X_\gamma)
\end{equation}
and thus any trace can be written as a symmetric trace plus lower degree traces contracted to structure constants. the previous section notes that we only get irreducible symmetric traces in certain dimensions.\\
Borrowing some of the notation from the previous section, if we write a symmetric trace with $k$ adjoint indices as a white vertex with $k$ thin lines coming out of it and a structure constant as a black vertex with $3$ thin lines coming out of it, the diagrams described above are equivalent to forests where each tree has one white vertex and some number of black vertices. By the discussion in the section on relations, the white vertices in turn must have degree $k = e_i+1$ for the exponents $e_i$ of $G$.\\
Note that for $B_r$ and $C_r$, the exponents are all odd, so the corresponding $k$ are all even. The white vertices can thus be replaced by traces over the adjoint representation rather than the defining representation $V$ using the relations described in [RSV], and the branches of black vertices can then be read as permutations of indices of traces in the adjoint representation, rather than permutations of indices of traces in the defining representation. Thus for $B_r$ and $C_r$, the adjoint tensor invariants can be expressed as tensor products of traces over the adjoint representation, with permutation of the indices. Because of the Killing form, the discussion in the section on $B_r$ also applies to the adjoint representation, so that traces of odd degree over the adjoint representation vanish identically. Hence for $A_r$ and $D_{2k+1}$, which both have even exponents, we cannot use the adjoint representation for everything, and we already noted that we need at least two representations for $D_{2k}$.
\section{Kirillov's Family Algebras and the Exceptional Cases}
In 2001 Kirillov defined a set of algebras that he calls family algebras. He considers a representation $V$ of $G$ and defines
\begin{equation}C_V(\mg) = (End(V)\otimes S(\mg))^G\end{equation}
The original intent was to compute the behavior of the $G$-harmonic polynomials of Kostant, applying algebraic methods in place of the combinatorial methods already known.\\
For the case where $G$ is classical, $V$ is a symmetric, antisymmetric, or tensor power of the adjoint representation, we can use the above theorem to describe the family algebra $C_V(\mg)$, noting that 
\begin{equation}C_V(\mg) = (V\otimes V^\vee \otimes S(\mg))^G \subset T(\mg)^G\end{equation}
Hence the generators of the tensor algebra describe the possible elements of the family algebra, which can then be described in terms of generators within the family algebra (using the family algebra's multiplication structure rather than tensor product). For instance, the author has used the main theorem to compute the structure of $C_{adj}(\mg)$ for $\mg$ classical and $adj$ the adjoint representation.\\
The author has also examined $C_{adj}(\mg)$ for $\mg$ exceptional. Just as the elements of the family algebras in the classical cases could be described as tensor products of traces, or alternatively as trees in the sense of the previous section, so could the elements of the family algebras in the exceptional cases. Thus the author conjectures that the main theorem may hold for the exceptional Lie algebras, with specific choices for the representation used to define $M$:
\begin{conjecture}
For $G_2$, let $(V,\pi)$ be the $7$-dimensional representation. For $F_4$, let $(V,\pi)$ be the $26$-dimensional representation. For $E_6$ let $(V,\pi)$ be the $27$-dimensional representation. For $E_7$ let $(V,\pi)$ be the $56$-dimensional representation. For $E_8$ let $(V,\pi)$ be the $248$-dimensional representation. Defining traces as above, $T(\mg)^G$ is generated by traces, allowing for permutation of the indices.
\end{conjecture}
Unfortunately, the argument used for the classical Lie algebras does not work here. The representations of the exceptional groups carry irreducible invariants of degree $3$ and $4$ that don't reduce or cancel. In particular, $G_2$'s and $E_8$'s defining representations carry antisymmetric cubic invariants, $F_4$'s and $E_6$'s defining representations carry symmetric cubic invariants, and $E_7$'s defining representation carries a symmetric quartic invariant. Thus, although all of the representations of a given exceptional group can be embedded into the tensor powers of the defining representation of that group, it is not clear that the invariants in the defining representation can be reduced to just Kronecker deltas.\\
Note that except for $E_6$, the exceptional Lie algebras have only odd exponents, and hence, again by [RSV], we can replace symmetric traces over $V$ with symmetric traces over the adjoint representation. Like $A_r$ and $D_r$, $E_6$ requires the use of at least one other representation.

\end{document}